\magnification\magstep1
\overfullrule = 0pt

\def\n{\noindent}
\magnification\magstep1
\overfullrule = 0pt

\def\n{\noindent}
\def\qed{{\hfill{\vrule height7pt width7pt
depth0pt}\par\bigskip}} 
\def\vp{\varepsilon}
\def\pf{\medskip\n {\bf Proof.}~~}

\def\n{\noindent}

\def\vp{\varepsilon}
\def\pf{\medskip\n {\bf Proof.}~~}

\def\CC{ \;Ê{}^{ {}_\vert }\!\!\!{\rm C}}

\def\CC{ \;Ê{}^{ {}_\vert }\!\!\!{\rm C}}

\def\CC{ \;Ê{}^{ {}_\vert }\!\!\!{\rm C}}

\def\NN{{\rm I}\!{\rm N}}

\def\vp{\varepsilon}
\def\pf{\medskip\n {\bf Proof.}~~}

\def\CC{ \;Ê{}^{ {}_\vert }\!\!\!{\rm C}}

\centerline{\bf Remarks on the similarity degree}
\centerline{\bf of an operator algebra}\bigskip
\centerline{ by Gilles Pisier\footnote*
{Supported in part by the NSF   and
 the Texas Advanced Research Program 010366-163.}}
\centerline{ Texas A\&M University}
\centerline{College Station, TX 77843, U. S. A.}
\centerline{and}
\centerline{Universit\'e Paris VI}
\centerline{Equipe d'Analyse, Case 186, 75252}
\centerline{ Paris Cedex 05, France}\bigskip
 
\centerline{ \bf  Abstract} \bigskip\bigskip
The ``similarity" degree of a unital operator 
 algebra $A$ was  defined and studied in two recent papers of ours, where
 in particular we showed that it coincides with the ``length" of an operator algebra. 
This paper  brings several complements: we give direct proofs (with slight improvements)
of several known facts on the length which were only known via the degree,
and we show that the length of a type $II_1$ factor with property $\Gamma$
is at most 5, improving on a previous bound ($\le 44$) due to E. Christensen.

\bigskip\bigskip

\centerline{********}

\vfill\eject
\bigskip\bigskip

The similarity degree of a unital operator 
 algebra $A$ is defined in [P1] as the 
smallest $\alpha\ge 0$ for which there is a constant $C$ such that any bounded 
morphism (= unital homomorphism) $u\colon \ A\to B(H)$ satisfies 
$$\|u\|_{cb} \le C\|u\|^\alpha.\leqno (1)$$
On the other hand, an operator algebra $A\subset B({\cal H})$ is said to be of 
length $\le d$ if there is a constant $K$  such that, for any $n$ and any $x$ 
in $M_n(A)$, there is an integer $N = N(n,x)$ and scalar matrices $\alpha_0 \in 
M_{n,N}(\CC)$, $\alpha_1\in M_N(\CC),\ldots, \alpha_{d-1}\in M_N(\CC)$, 
$\alpha_d \in M_{N,n}(\CC)$ together with diagonal matrices $D_1,\ldots, D_d$ in 
$M_N(A)$ satisfying
$$\left\{\eqalign{&x = \alpha_0 D_1 \alpha_1 D_2 \ldots D_d\alpha_d\cr
&\prod^d_0 \|\alpha_i\| \prod^d_1 \|D_i\| \le K\|x\|.}\right.\leqno (2)$$
We denote by $\ell(A)$ the smallest $d$ for which this holds and we call it the 
``length'' of $A$ (so that $A$ has length $\le d$ is indeed the same as $\ell(A) 
\le d$). It is easy to see that if $\ell(A) \le d$ as in (2) above then for any 
bounded homomorphism $u\colon \ A\to B(H)$ we have $\|u\|_{cb} \le K\|u\|^d$. 
In [P1], we proved the following basic result in the converse direction.

\proclaim Theorem 1. For any unital operator
 algebra $A$, we have $$d(A) = \ell(A).$$

\n {\bf Remark 2.} Let $\pi\colon \ A\to B(H)$ be a completely contractive 
homomorphism and let $\delta\colon \ A\to B(H)$ be a $\pi$-derivation (i.e.\ we 
have $\delta(xy) = \delta(x)  \pi(y)  + \pi(x) \delta(y)$ for all $x,y$ in $A$). 
Then it is easy to see that if $\ell(A)\le d$  as in (2) above then 
$\|\delta\|_{cb} \le K d\|\delta\|$. The following result is a converse to this 
essentially due to Kirchberg [K], with an improvement observed in [P1].

\proclaim Proposition 3. Let $A$ be a unital $C^*$-algebra. Assume there is a 
constant $C$ such that for any $\pi$ and $\delta$ we have $\|\delta\|_{cb} \le 
C\|\delta\|$. Then $d(A) \le C$ (hence $d(A)$ is at most equal to the integral 
part of $C$). 

There is a lot of information known on the value
of $d(A)$ for various examples  due to the works
of E.~Christensen and Uffe Haagerup ([C1-4, H]).
In this note, we try to  give direct and explicit
factorizations for the corresponding results for
$\ell(A)$. 
Because of Theorem~1, our results are essentially
already known but they shed  some light on the
meaning of
$\ell(A)$, and we obtain some  slight
improvements.\medskip

\n {\bf Remark 4.} (i) \ Fix an integer $n$ and consider $x$ in $M_n(A)$. We 
denote 

$$ \|x\|_{(d)}=\inf\{
 \prod\limits^d_{i=0} \|\alpha_i\| 
\prod\limits^d_{i=1} \|D_i\|\}$$
where the infimum runs  over all possible
representations   $x = \alpha_0 
D_1\alpha_1\ldots D_d\alpha_d$ of the form
appearing in (1) above. Clearly 
$$\|x\|_{M_n(A)} \le \|x\|_{(d)}.$$ 
Sometimes to avoid possible confusion in the sequel we will denote
$\|x\|_{(d)}$ instead by $\|x\|_{(d,A)}$.

\n (ii) \ It is rather easy to check, say if $A$ is unital, that $\|~~~\|_{(d)}$ 
is an equivalent norm on $M_n(A)$ for each fixed $n$. Actually, for any $x$ in 
$M_n(A)$ we have
$$\|x\|_{(1)} \le n \|x\|,$$
but we will content ourselves with
$$\|x\|_{(1)} \le n^2 \sup_{ij} \|x_{ij}\|\leqno (3)$$
which is easy to show: \ we take $N=n^2$, we introduce a bijection $\sigma\colon 
\ [1,2,\ldots, n]^2 \to [1,\ldots, N]$ and we let $D$ in $M_N(A)$ be the 
diagonal matrix with entries $D_k = x_{\sigma^{-1}(k)}$. Then we set $\alpha_0 = 
\sum\limits^n_{i,k=1} e_{i\sigma(i,k)}$ and $\alpha_1 = \sum\limits^n_{k,j=1} 
e_{\sigma(k,j)j}$. Clearly we have $\alpha_0D\alpha_1 = x$, $\|\alpha_0\| \le 
n$, $\|\alpha_1\|\le n$ and $\|D\| = \sup\limits_{ij} \|x_{ij}\|$ so that (3) 
follows. Note that if $A$ has a bounded approximate unit (say on the right) (3) 
implies an estimate $\|x\|_{(d)} \le C_n \sup\limits_{ij} \|x_{ij}\|$ for some 
constant $C_n$ and $\|~~~\|_{(d)}$ are equivalent norms on $M_n(A)$.

\n (iii) \ When $A$ is unital (or has a contractive approximate unit, which is 
the case whenever $A$ is a $C^*$-algebra) we have
$$\|x\|_{(d+1)} \le \|x\|_{(d)}.\leqno \forall x \in M_n(A)$$
Moreover, if $A$ is a $C^*$-algebra, the results of [BP] show that
$$\|x\|_{M_n(A)} = \lim_{d\to \infty} \downarrow \|x\|_{(d)}.\leqno \forall x\in 
M_n(A)$$

\proclaim Lemma 5. Let $A\subset B(H)$ be a closed subalgebra of $B(H)$. We 
assume that there are elements $p,q$, $a_i, b_i,c_i,d_i$ $(1\le i\le n)$ in $A$ 
such that
$$\eqalign{a_ib_j = \delta_{ij}p,\quad c_id_j =
\delta_{ij}q ,\quad\left\|\sum
b_ib^*_i\right\|\le 1 ,\quad
\left\|\sum c^*_jc_j\right\|\le 1 
.}$$
Let then $W_1,W_2$ be two unitary $n\times n$ matrices such that
$$|W_1(i,j)| = |W_2(i,j)| = n^{-1/2}.\leqno \forall i,j =1,\ldots, n$$
Then for any $x = [x_{ij}]$ in $M_n(A)$ we have the following factorization
$$[px_{ij}q] = D_1W_1 D_2W_2D_3$$
where $D_1,D_2,D_3$ are diagonal
 matrices with entries in $A$ such that  
$$ {\|D_1\| 
= \sup\|a_i\|,\quad \|D_3\|= \sup\|d_i\|,\quad
\hbox{and}\quad \|D_2\|  \le \|x\|_{M_n(A)}.}$$

\pf Indeed, let
$$\vp^1_{ij} = n^{1/2} \overline{W_1(i,j)} \quad \hbox{and}\quad \vp^2_{ij} = 
n^{1/2} \overline{W_2(i,j)}.$$
We set $D_1(i) = a_i, D_3(i) = d_i$ and
$$D_2(k) = \sum_{ij} \vp^1_{ik} b_i x_{ij} c_j \vp^2_{kj}.$$
Then it is easy to check that for any $\xi,\eta$ in the unit  ball of $H$ we 
have
$$|\langle D_2(k) \xi,\eta\rangle| \le \|x\| \left(\sum 
\|c_j\xi\|^2\right)^{1/2} \left(\sum_i \|b^*_i\eta\|^2\right)^{1/2} \le \|x\|.$$
Hence $\|D_2(k)\| \le \|x\|$.
Moreover a simple calculation shows that
$$\eqalign{(D_1W_1D_2W_2D_3)_{ij} &= a_i \sum\nolimits_k W_1(i,k) D_2(k) 
W_2(k,j)d_j\cr
&= \sum_k W_1(i,k) [a_iD_2(k)d_j] W_2(k,j)\cr
&= \sum_k W_1(i,k) [\vp^1_{ik} px_{ij}q\vp^2_{kj}] W_2(k,j)\cr
&= px_{ij}q \sum_k n |W_1(i,k)|^2 |W_2(k,j)|^2\cr
&= px_{ij}q.}$$
\qed
\medskip

\proclaim Corollary 6.  Let $W_1,W_2$ be as in Lemma~5 and let $A = B(H)$ with 
$\dim H = \infty$. Then for any $x$ in the unit ball of $M_n(A)$, there are   diagonal matrices
$D_1,D_2,D_3$ in the unit ball of 
$M_n(A)$ such that
$$x = D_1 W_1 D_2 W_2 D_3.$$
More generally, the same is true for any $C^*$-algebra $A$ which contains, for 
each $n$, a set of isometries $s_1,\ldots, s_n$ such that $s^*_is_j = 
\delta_{ij}1$.

\pf We simply take $a_i = s^*_i$, $b_i=s_i$, $c_j=s^*_j$ and $d_i=s_i$ in 
Lemma~5.
\qed
\medskip

\n {\bf Remark 7.} For any $C^*$-algebra, we have
$$\ell(A) \le \ell(A^{**}).$$
Indeed, let $\pi$ be the universal representation
of $A$, and let $M=\pi(A)''$ so that $A^{**}=M$.
 Fix $x$ in the unit ball of $M_n(A)$. Let $\hat x \in M_n(M)$ be 
the same element viewed in $M_n(M)$. If $\ell(M)\le d$, we can write 
$\hat x = \alpha_0 D_1\alpha_1D_2\ldots D_d\alpha_d$ with $\|D_i\| \le 1$,
$D_i$ diagonal with entries in $M$,  
$\prod 
\|\alpha_i\|\le K$. Then let $(D^\alpha_i)$ be a
net in the unit ball of $A$ tending  in the strong
operator topology (in short sot) to $D_i$. Then
if $x^\alpha = 
\alpha_0 D^\alpha_1 \alpha_1\ldots D^\alpha_d\alpha_d$ we have $x^\alpha_{ij}\to 
x_{ij}$ in sot, hence (since we are dealing
with the universal representation) in the weak topology of $A$, hence after passing to the
convex hull
$x^\alpha_{ij} \to  x_{ij}$ in norm, which implies by Remark~4.ii that $\|x\|_{(d)} \le 
\lim\limits_\alpha \|x^\alpha\|_{(d)} \le K$.\qed

\proclaim Theorem 8. Let $A$ be an operator algebra. We have $\ell(A) \le 3$ in 
the following cases:
\item{(i)} $A = B(H)$,
\item{(ii)} $A$ is a  $C^*$-algebra without
tracial states,
\item{(iii)} $A = K \otimes_{\rm min} B$ where $B$ is an 
arbitrary unital operator algebra.
\medskip

\pf (i) follows  from Corollary 6. To check (ii), note that if $A$ has no 
tracial states then it is well known that (for any $n\ge 1$) $A^{**}$ contains 
isometries $s_1,\ldots, s_n$ such that $s^*_is_j = \delta_{ij}=1$. Hence by 
Corollary~6, $\ell(A^{**}) \le 3$ and by Remark~7 we obtain $\ell(A) \le 3$. To 
check (iii), note that by (3) it suffices to be able to factorize all elements 
$x$ 
of a dense subset of the unit ball of $M_n(A)$. Hence we may assume that the 
entries $x_{ij}$ of $x$ lie in a dense linear subspace $V\subset A$. We will use 
the linear subspace $V\subset K \otimes_{\rm  min} B$ spanned by 
$[e_{ij}\otimes b]$. Then consider $x$ in the unit ball of $M_n(V)$, so that 
there is an integer $m$ such that if $p,q$ in $V$ are defined as $p=q  = 
\sum\limits^m_{i=1} e_{ii} \otimes 1_B$, we have $x_{ij} = px_{ij}q$. Thus, if 
we set $a_i = ps^*_i \otimes 1_B$, $d_i = s_ip \otimes 1_B$, $b_i = s_i \otimes 
1_B$ and $c_i = s^*_i \otimes 1_B$, we obtain $\|x\|_{(3)} \le \|x\|_{M_n(A)}$ 
when $A = K\otimes_{\rm min} B$ and we conclude that $\ell(K \otimes_{\rm min} 
B) \le 3$.\qed

\proclaim Theorem 9. Let $A\subset B(H)$ be a $C^*$-subalgebra, generating a 
von~Neumann algebra $M$. Assume that $M$ is a $II_1$-factor. Then $\ell(M) \le 
\max\{\ell(A), 3\}$.

\n {\bf Corollary 10.} Let $R$ be the hyperfinite $II_1$-factor. Then $\ell(R) 
\le 
3$. 

The last two statements are easy consequences of
the following.

\proclaim Lemma 11. Let $x = (x_{ij})$ in the unit ball of $M_n(M)$ be such that 
$\sum\limits_{ij} \|x_{ij}\|^2_2 <\vp^2$, then
there are projections $p,q$ in 
$M$ with $\tau(p) \le {1\over n}$, $\tau(q) \le {1\over n}$ such that
$$x_{ij}  = px_{ij}q + y_{ij}$$
with $\|y_{ij}\|\le 2\vp\sqrt n$.

\pf Consider $a = \left(\sum\limits_{ij} x^*_{ij}x_{ij}\right)^{1/2}$. Then if 
$E_a$ 
the spectral measure of $a$, we let $q = E_a [\vp \sqrt n, \infty]$ so that 
$\tau(q) \le {1\over n}$
$$x_{ij} = x_{ij}q + x_{ij}(1-q)$$
so that
$$\eqalign{\|x_{ij}(1-q)\|^2 &= \|(1-q)x^*_{ij} x_{ij}(1-q)\|\cr
&\le \|(1-q) a^2(1-q)\|\cr
&\le \vp^2n.}$$
Hence $\|x_{ij}(1-q)\| \le \vp \sqrt n$. Now let $b = \left(\sum\limits_{ij} 
x_{ij}x^*_{ij}\right)^{1/2}$ and $p = E_b[\vp \sqrt n,\infty]$. We have
$$x_{ij}q = px_{ij}q + (1-p)x_{ij}q$$
and again $\tau(p) \le {1\over n}$ but
$$\|(1-p)x_{ij}q\| \le \|(1-p)x_{ij}\| \le \vp \sqrt n.$$
\qed

\n {\bf Remark.} As a consequence of Lemma 11,
 if $C(d)$ denotes the unit ball of 
the norm $\|~~~\|_{(d)}$ on $M_n(M)$ (see
Remark~4) and if $d\ge 3$ we have: \ 
$$\overline{C(d)}^{L_2(M_n(M))} \subset
3C(d).$$\medskip

\n {\bf Proof of Theorem 9.} By the Kaplansky density theorem, for any $\vp>0$, 
any $z$ in the unit ball of $M_n(M)$ can be  written as $z = z'+x$ with $z'$ in the 
unit ball of $M_n(A)$ and $\sum \|x_{ij}\|^2_2 < \vp^2$. With the notation of 
Remark~5, if $d = \ell(A)$ we have $\|z'\|_{(d)} \le K$ for some fixed constant 
$K$ (independent of $n$). Thus we are reduced to the factorization of $x$. But 
with the notation of Lemma~11 we have $x_{ij} = px_{ij}q + y_{ij}$ and Lemma~4 
ensures that $\|(px_{ij}q)\|_{(3)} \le 2$. Thus we are reduced to estimate $y = 
(y_{ij})$, but by (3) we have
$\|y\|_{(1)} \le 2\vp n^{5/2}$, hence we finally conclude that if $d\ge 3$ we 
have 
$$\|z\|_{(d)} \le K + 2 + 3\vp n^{5/2}$$
and if $d<3$ we have the same majorization for $\|z\|_{(3)}$. Thus we obtain 
Theorem~9.\qed \medskip

\n {\bf Remark 12.}  
It is proved in [P1] that $d(A) 
\le 2$ implies that $A$ is ``semi-nuclear" in the following sense:
for any $*$-representation
$\pi: A \to B(H)$, the generated von Neumann algebra $M=\pi(A)''$ is injective whenever it is is   semi-finite.
Note that nuclear implies semi-nuclear.
By well known results, if $G$ is any discrete group, and if either $A= C^*(G)$  
   or  $A=C_\lambda^*(G)$, then $A$ is semi-nuclear iff $G$ is amenable,
or equivalently  iff $A$ is nuclear. In general it seems to be open whether
conversely semi-nuclear
implies nuclear. However, the results of [A] imply
that $B(H)$ is not semi-nuclear (here $\dim(H)=\infty$).
As pointed out to me by Narutaka Ozawa, it is easy to adapt the argument
in [A] to show that the hyperfinite $II_1$ factor $R$  is not semi-nuclear,
and  actually that no  
$II_1$ factor can be semi-nuclear. Thus
we have $$\ell(M)\ge 3$$
for any $II_1$ factor $M$.

 Let $M$ be a $II_1$-factor with (Murray and
von~Neumann's) property 
$\Gamma$. This means that there is a 
net of unitaries $(u_\alpha)$ in $M$  with 
zero trace which are asymptotically central, i.e.\ are such that $\|u_\alpha t - 
tu_\alpha\|_2 \to 0$ for any $t$ in $M$.

By a result of Dixmier [D], we can then find ``many'' asymptotically central 
elements, in particular for any $n$ there exists a net $(p^\alpha_1,\ldots, 
p^\alpha_n)$ of orthogonal decompositions of the identity in $M$ with 
$\tau(p^\alpha_i) = {1\over n}$ for all $i,\alpha$ and such that 
$(p^\alpha_i)_\alpha$ is asymptotically central for each $i=1,\ldots, n$. In 
particular, we have for any $t$ in $M$, $\lim\limits_\alpha \left\|t - \sum^n_1 
p^\alpha_i tp^\alpha_i\right\|_2  = 0$
 (indeed note that $\|t - [ptp + (1-p)
t(1-p)]\|_2  =  \|t -[ptp + t-pt - tp +
ptp]\|_2 =  \|p(tp-pt) + (pt-tp)p\|_2). $)
Now fix $x$ in the unit ball of $M_n(M)$ and let $x^\alpha_{ij} = 
\sum\limits^n_{m=1} p^\alpha_m x_{ij} p^\alpha_m$. We have 
$\|x^\alpha\|_{M_n(M)} \le \|x\|_{M_n(M)}$ and
$$\lim_{\alpha\to \infty} \|x-x^\alpha\|_{L^2(M_n(M))} =0.$$
This allows us to improve 
the main result of [C4], as follows (the estimate given
in [C4] is $d(M)\le 44$).

\proclaim Theorem 13. If $M$ is $II_1$-factor with property $\Gamma$ then 
$\ell(M) \le 5$.

\pf By the preceding remarks and by Lemma~11, it suffices to show that if 
$p_1,\ldots, p_n$ are disjoint projections in $M$ with $\tau(p_i) = {1\over n}$ 
we have
$$\left\|\left[ \sum^n_{m=1} p_m x_{ij} p_m \right]_{ij}\right\|_{(5)} \le 
\|x\|_{M_n(M)}.$$
This is an immediate consequence of the next
lemma, where we use freely the notation
introduced in Remark 4.

\proclaim Lemma 14. Let $X_1,\ldots, X_n$ be in $M_n(M)$ with $\|X_m\|_{(d)} \le 
1$ for any $m=1,\ldots, n$. Then if we let
$$y_{ij} = \sum^n_{m=1} p_m X_m(i,j) p_m$$
with $p_m$ as above, we have
$$\|y\|_{(d+2)} \le 1.$$

\pf We have
$$y = \alpha X \alpha^*$$
where $\alpha \in M_{n,n^2}(\CC)$ and $X\in M_{n^2}(M)$ are defined as:
$$\alpha = \sum \bar e_{1m} \otimes 1 \otimes p_m$$
and $X = \sum e_{mm}\otimes X_m$. Here $\bar e_{1m}$ denotes the canonical basis 
of $M_{1,n}(\CC)$ and we use the 
usual identifications $$M_{1,n} \otimes M_n = 
M_{n,n^2}\quad  \hbox{and} \quad M_n\otimes M_n =
M_{n^2}.$$

\n  It is easy to see, that in $M_{n^2}(M)$ we
have
$\|X\|_{(d)} \le \sup\limits_m 
\|X_m\|_{(d)}$. Thus it suffices to show that $\|\alpha\|_{(1)}\le 1$, since 
$\|\alpha^*\|_{(1)} \le 1$ follows by transposition. But the latter follows from 
the following identity
 $$\alpha = \sum^n_1 \bar e_{1m} \otimes 1 \otimes
p_m= \alpha_0 DW$$
 where $\alpha_0 = {1\over \sqrt n} \sum^n_1 \bar e_{1m} \otimes 1 \otimes 1$
$$\eqalignno{D &= \sum\nolimits^n_1 e_{mm} \otimes 1 \otimes \sum^n_{i=1} 
\overline W_{im} \sqrt n \ p_i\cr
\noalign{\hbox{and}}
W &= \sum_{ij} W_{ij} e_{ij} \otimes 1 \otimes
1.}$$ Indeed, we have 
$$\eqalign{\alpha_0 DW &= {1\over \sqrt n}
 \sum^n_{m=1} \sum^n_{j=1} \bar e_{1j} 
\otimes 1 \otimes \sum^n_{i=1} \sqrt n \
\overline W_{im} W_{mj} p_i\cr &= \sum_j
\bar e_{1j} \otimes 1 \otimes \sum_i (W^*W)_{ij}
p_i\cr &= \sum^n_{j=1} \bar e_{1j}  \otimes 1
\otimes p_j = \alpha.}$$ Whence we obtain
$$\|\alpha\|_{(1)} \le \|\alpha_0\| \ \|D\| \ \|W\| \le 1.$$
\qed

\vfill\eject
We will now discuss the following
 which goes back to 1955. \medskip

\n {\bf Kadison's conjecture} ([K]). Any
$C^*$-algebra has the following similarity
property (SP). Every bounded homomorphism
$u\colon \ A\to B(H)$ is similar to a $*$-homomorphism, i.e.\ there is an
invertible $\xi\colon \ H\to H$ such that the homomorphism $u_\xi$ defined by
$u_\xi(x) = \xi^{-1} u(x)\xi$ $(\forall x\in A)$
is a
$*$-homomorphism, which means that $u_\xi(x^*) =
u_\xi(x)^*$ $(\forall x\in A)$.\medskip

 By a result due to Haagerup [H], it is
known that (if $u$ is unital) $u$ is similar to a $*$-homomorphism iff $u$ is
c.b. Moreover, we have
$$\|u\|_{cb} = \inf\{\|\xi^{-1}\| \ \|\xi\|\}$$
where the infimum runs over all $\xi$ for which $u_\xi$ is a $*$-homomorphism
(or equivalently for which $\|u_\xi\|  = 1$).

By the results of [P1], we have

\proclaim Proposition 15. If Kadison's conjecture is true than there is a
fixed $d_0$ such that any $C^*$-algebra has length $\le d_0$.

\pf By [P1] a $C^*$-algebra $A$ satisfies (SP) iff $\ell(A) < \infty$. Assume
that there are $C^*$-algebras $A_n$ such that $\ell(A_n) \to \infty$ when
$n\to \infty$. Then let $A =  \bigoplus\limits_n A_n$ be (say) the
$C^*$-algebra formed of sequences $x = (x_n)_n$ with $x_n\in A_n$ such that
$x_n\to 0$ when $n\to \infty$. Clearly $\ell(A)\ge\ell(A_n)$ for all $n$ hence
$\ell(A) = \infty$. Thus if Kadison's conjecture is correct, $\ell(A) <\infty$
for any $A$, whence Proposition~15 follows.\qed

Let $A$ be an operator algebra. For any set
$I$, we denote by $\ell_\infty(I,A)$ the
operator algebra formed of all bounded
families $(x_i)_{i\in I}$ with $x_i\in A$ for
all $i$ in $I$, equipped with the norm $\|x\|
= \sup\limits_{i\in I} \|x_i\|$. In the next
result, we show that the length of
$\ell_\infty(I,A)$ (with $I$ infinite) is
closely related to the possibility of
obtaining the factorization described in (2)
above, with the size $N$ and the scalar
factors depending only on $n$ and not of
$x\in M_n(A)$.

\proclaim Proposition 16. Let $A$ be an operator algebra. For any set $I$, we
denote by $\ell_\infty(I,A)$ the operator algebra formed of all bounded
families $(x_i)_{i\in I}$ with $x_i\in A$ for all $i$ in $I$, equipped with
the norm $\|x\| = \sup\limits_{i\in I} \|x_i\|$. Fix an integer $d\ge 1$. The
following assertions are equivalent.
\item{(i)} For any set $I, \ell(\ell_\infty(I,A)) \le d$.
\item{(ii)} \ For any countable set $I$,
 $\ell(\ell_\infty(I,A))\le d$.
\item{(iii)} There is a constant $K$ such that for
any $n$ there is an integer
$N=N(n)$ and scalar matrices of norm 1
$$\alpha_0\in M_{n,N}(\CC), \alpha_1\in M_N(\CC),\ldots, \alpha_{d-1} \in
M_N(\CC),  \alpha_d\in M_{N,n}(\CC)$$
 such that for any $x$ in $ M_n(A)$ there are
diagonal matrices $D_1,\ldots, D_d$ in $M_N(A)$ with $\prod\limits^d_1 \|D_i\|
\le K\|x\|$ and satisfying
$$x = \alpha_0D_1\alpha_1D_2\ldots D_d\alpha_d.$$

\pf Note that we have canonically
$$M_n(\ell_\infty(I,A)) = \ell_\infty(I,M_n(A)).$$
Assume (iii). Consider $x$ in
$M_n(\ell_\infty(I,A))$. Let $(x(i))_{i\in I}$ be
the associated element in $\ell_\infty(I,M_n(A))$
with 
$$\|x\| = \sup_{i\in I} \|x(i)\|.$$
By (iii) we can find for each $i$ in $I$ diagonal
elements $D_1(i),\ldots, D_d(i)$ in $M_N(A)$ such
that $x(i)  = \alpha_0D_1(i)\ldots
D_d(i)\alpha_d$. Let $D_1,\ldots, D_d$ be the
corresponding (diagonal) elements of
$M_N(\ell_\infty(I,A))$. Then we clearly
 have $  x = \alpha_0 D_1\alpha_1
\ldots D_d\alpha_d$ and (2) holds. So we obtain
(i).

Conversely, assume (i). Let $I$ be the unit ball of $M_n(A)$. Let $x\colon \
I\to M_n(A)$ $(i\to x(i)$) be the inclusion mapping. Clearly $x\in
\ell_\infty(I, M_n(A))$ with $\|x\| = 1$. We can also view $x$ as an element
of $M_n(\ell_\infty(I,A))$. Then, if (i) holds we can find $\alpha_0,\ldots,
\alpha_d$ with norm 1 and $D_1,\ldots, D_d$ diagonal in
$M_N(\ell_\infty(I;A))$ such that $x = \alpha_0D_1\ldots D_d\alpha_d$ and
$\Pi\|D_i\| \le K\|x\| \le~K$. (Here $N = N(n,x)$ but since $x$ is 
a fixed canonical element, actually $N$ depends
only on $n$.) Taking the $i$-th coordinate, we
obtain
$$x(i) = \alpha_0D_0(i)\alpha_1,\ldots D_d(i)\alpha_d.\leqno \forall i\in I$$
Hence we conclude (by homogeneity) that (iii)
holds.

\n This shows
that (i)~$\Leftrightarrow$~(iii). Since (i)~$\Rightarrow$~(ii) is trivial, it
remains only to show that (ii)~$\Rightarrow$~(iii). Assume (ii). Fix an
integer $n$. Then we may observe that $\ell_\infty(\NN,A)$ is of length $\le
d$ with $K$ (independently of $n$ or $x$) and $N(n,x)$ as defined before (2)
but moreover with matrices $\alpha_0,\alpha_1,\ldots, \alpha_d$ with {\it
rational\/} coefficients.
Indeed, a simple density argument 
(it is convenient here to  invoke (3)) establishes
this fact (perhaps at the cost of a $K$ slightly
worse than the original one).

\n Then let $S_N$ be the set of all $(d+1)$-tuples
$(\alpha_0,\ldots, \alpha_d)$ with rational
coefficients in $M_{n,N}\times M_{N,N}\times
\cdots
\times M_{N,n}$ with
$\|\alpha_0\|,\ldots, \|\alpha_d\|\le 1$, and let
$$I = \bigcup_{N\ge 1} \{N\}\times S_N.$$
Clearly, $I$ is countable.

We claim that there is a $p$ in $I$ say $p = (N, (\alpha_0,\ldots, \alpha_d))$
with $(\alpha_0,\ldots, \alpha_d)\in S_N$ such that (iii) holds relative to
$N$ and $(\alpha_0,\ldots,\alpha_d)$ (i.e.\ the same $N$ and the same
$(\alpha_0,\ldots, \alpha_d)$ work for any $x$ in the unit ball of $M_n(A)$).

\n Indeed, if we assume otherwise. Then for any
$p$ in $I$, there is $x_p$ in the unit ball of
$M_n(A)$ such that whenever we have an equality
$$x_p = \alpha_0 D_1\alpha_1\ldots D_d\alpha_d$$
with $D_1,\ldots, D_d$ as in (iii) then we must have
$$\prod^d_{j=1} \|D_j\| > K.\leqno (4)$$
Let $x = (x_p)_{p\in I}$. Note that $x$ 
is in the unit ball of
$\ell_\infty(I; M_n(A)) =
M_n(\ell_\infty(I,A))$. Applying our original
observation about {\it rational\/} coefficients,
we find that there is an $N$ and $q =
(\alpha_0,\ldots, \alpha_d)$ in $S_N$ such that
$x$ can be written as
$x = \alpha_0D_1\ldots D_d\alpha_d$ with
 $D_1,\ldots, D_d$ diagonal in
$M_N(\ell_\infty(I,A))$ such that
$\prod\|D_j\|\le K$. In particular, if we
restrict this equality to the $p$-th coordinate
of $x$ with $p=(N,q)$ we find a factorization of
the form $x_p  = \alpha_0D_1(p)\ldots
D_d(p)\alpha_d$ which contradicts (4). 

This proves the above claim, and thus concludes the proof
that (ii)~$\Rightarrow$~(iii).\qed

Taking into account Proposition~16, a close look at the proof of Theorem~9
immediately yields:

\proclaim Corollary 17. Let $M$ be a
$II_1$-factor with property $\Gamma$. Then we have
$\ell(\ell_\infty(I,M)) \le 5$ for any set $ I$.

Let $B\subset A$ be an inclusion between operator algebras. In order to
compare the norms $\|~~\|_{(d)}$ defined above in Remark~4 for $A$ and for
$B$, we denote the respective norms by $\|\cdot\|_{(d,A)}$ and
$\|~~\|_{(d,B)}$.

\proclaim Sublemma 18. Let $B$ 
be an operator algebra and let $A =
M_n(B)\approx B\otimes M_n$. Fix $r,s$ with $1\le
r,s\le n$. Let $j_{r,s}\colon \ B\to A$ be the
(``inclusion'') mapping defined by $j_{r,s}(b) =
b\otimes e_{rs} (b\in B)$. Then, for any $x$ in
$M_n(B)$ we have ($r,s$ remaining fixed) 
$$\|[j_{r,s}(x_{ij})]\|_{(3,A)} \le \|x\|_{M_n(B)}.$$

\pf We apply Lemma~5 with $q = p = 1\otimes
e_{rs}$ $a_i = 1 \otimes e_{ri}$,
$b_j = 1 \otimes e_{js}$, $c_i = 1 \otimes
e_{ri}$, $d_j = 1 \otimes e_{js}$.\qed

\proclaim Sublemma 19. Again let $A =
M_n(B)\approx  B\otimes M_n$ as before and let
$j\colon \ B\to A$ be the mapping defined by
$J(b) = b\otimes 1$. Then for any $x = (x_{ij})$
in
$M_n(B)$ we have
$$\|[J(x_{ij})]\|_{(5,A)} \le \|x\|_{M_n(B)}.$$

\pf This is the same argument as for Lemma~14. We can write
$$J(x_{ij}) = x_{ij} \otimes 1 = \sum^n_{m=1} p_m
X_m(i,j)p_m$$ where $p_m = 1\otimes e_{mm} \in A$
and
$$X_m(i,j) = x_{ij} \otimes e_{mm}  =
J_{mm}(x_{ij}).$$ Then arguing as for Lemma~14,
we obtain Sublemma~18.\qed\medskip

\n {\bf Remark 20.} Note that the factors $W_1,W_2,W$ (and their sizes) which
appear when we spell out explicitly the factorizations corresponding to
Sublemmas~18 and 19 depend only on $n$ and not on $x$. 

Hence we obtain:

\proclaim Proposition 21. Let $C$ be the CAR algebra $C = (M_2)^{\otimes \NN}$
or any infinite $C^*$-tensor product of matrix
algebras. Then for any set $I$, we have
$\ell(\ell_\infty(I,C)) \le 5$.

\pf Consider $x$ in the open unit ball of $M_n(C)$. By density we may assume
that all entries $x_{ij}$ belong to $C_N\otimes 1 \simeq C_N$ where $C_N = M_2
\otimes\cdots\otimes M_2$ ($N$-times). Now assume without loss of generality
that $n = 2^k$. Note that the inclusion $C_N\to C$ can be factored as $C_N
{\buildrel J\over\longrightarrow } C_N \otimes M_{2^k} {\buildrel
\pi\over\longrightarrow } C$ where $J(b) = b\otimes 1$ as above. Thus since
$\pi$ is a $*$-homomorphism we have $\|x\|_{(5,C)} \le \|x\|_{M_n(C)}$. By
Remark~20 and Proposition~16, we conclude that $\ell(\ell_\infty(I,C)) \le
5$.\qed

The $II_1$ factor associated with the free group with at least two generators
is a typical example of one failing property $\Gamma$, and it might be
a counterexample to Kadison's conjecture. But actually, we feel
that the following should be true.\medskip

\n {\bf Conjecture.} Let $M$ be the von Neumann algebra formed
of all norm-bounded  sequences $(x_n)$ with $x_n\in M_n$ for each $n$,
equipped with the sup-norm and let $N=\ell_\infty(\NN,M)$.
Then $N$ (and perhaps even $M$) is a counterexample to   Kadison's conjecture.
In other words, its ``length" is infinite.

\bigskip\bigskip
\vfill\eject
\centerline {\bf References}
\bigskip
 \item{[A]}
  J. Anderson. Extreme points in sets of positive linear maps on
${\cal B}({\cal H})$. J. Funct. Anal. 31 (1979) 195-217. 
 
 \item{[BP]} D. Blecher and V. Paulsen. 
 Explicit construction of universal operator
 algebras and applications to polynomial
 factorization. Proc. Amer. Math. Soc.   112
 (1991) 839-850.

\item{[C1]}    E. Christensen.   Extensions of derivations. 	J.
Funct. Anal.  27 (1978) 234-247.

\item{[C2]}    E. Christensen.   Extensions of derivations II.
Math. Scand. 50 (1982)  111-122.

\item{[C3]}   E. Christensen. On non self adjoint
representations of operator algebras  Amer. J. Math.
103 (1981) 817-834.

\item{[C4]}   E. Christensen. Similarities of $II_{1}$
factors with property $\Gamma$.  Journal Operator Theory 
15 (1986) 281-288.

 \item{[D]} J. Dixmier.   Quelques propri\'et\'es
des suites centrales dans les facteurs de type
${\rm II}\sb{1}$. (French) Invent. Math. 7 (1969)
215-225. 

\item{[H]}  U. Haagerup.  Solution of the
similarity problem for cyclic representations of
$C^*$-algebras.  Annals of Math.  118 (1983),
215-240.

\item{[Ka]} R.  Kadison  On the orthogonalization
of operator representations. Amer. J. Math.  77
(1955) 600-620.

\item{[Ki]} E. Kirchberg.  The derivation and the similarity problem 
are equivalent.
 J. Operator Th. 36 (1996) 59-62.

\item{[Pa1]} V. Paulsen.   Completely bounded maps and
dilations.  Pitman Research Notes in
Math. 146, Longman, Wiley, New York, 1986.

\item{[P1]}  G. Pisier.  The similarity degree of
an operator algebra.
   St. Petersburg Math. J. 10 (1999) 103-146.
 
\item{[P2]}   $\underline{\hskip1.5in}$. Joint
similarity problems and the generation of
operator algebras with bounded length. Integr.
Equ. Op. Th.
 31 (1998) 353-370.

\item{[P3]}   $\underline{\hskip1.5in}$.
Similarity problems and completely bounded maps.
Springer Lecture notes 1618 (1995).

 \bye

\n{\bf Remark.} Instead of proving this for
$N$, we feel it might be easier to prove    
for $N_\infty=\ell_\infty(\NN,N)$.

\vskip12pt

Texas A\&M  University

College Station, TX 77843, U. S. A.

and

Universit\'e Paris VI

Equipe d'Analyse, Case 186,
 
75252 Paris Cedex 05, France
\end

\bye